\documentclass{birkjour}
 \newtheorem{thm}{Theorem}[section]
 
 \newtheorem{lem}[thm]{Lemma}
 
 \theoremstyle{definition}
 
 \theoremstyle{remark}
 \newtheorem{rem}[thm]{Remark}
 \newtheorem*{ex}{Example}
 \numberwithin{equation}{section}

\begin{document}

%
%
%
%
%
%
%
%
%

\title[Harmonic Univalent Functions Involving Pascal Distribution Series]
 {Some Connections Between Various Subclasses of Harmonic Univalent Functions Involving Pascal Distribution Series}

\author[Serkan \c{C}akmak]{Serkan \c{C}akmak}

\address{%
Department of Mathematics, Faculty of Arts and Science\\
Bursa Uludag University\\
16059 Bursa\\
Turkey}

\email{serkan.cakmak64@gmail.com}

\author{Elif Ya\c{s}ar}
\address{Department of Mathematics, Faculty of Arts and Science\br
Bursa Uludag University\br
16059 Bursa\br
Turkey}
\email{elifyasar@uludag.edu.tr}

\author{Sibel Yal\c{c}\i n}
\address{Department of Mathematics, Faculty of Arts and Science\br
Bursa Uludag University\br
16059 Bursa\br
Turkey}
\email{syalcin@uludag.edu.tr}

\author{\c{S}ahsene Alt\i nkaya}
\address{Department of Mathematics, Faculty of Arts and Science\br
Bursa Uludag University\br
16059 Bursa\br
Turkey}
\email{sahsene@uludag.edu.tr}



\date{January 1, 2004}

\begin{abstract}
In the present paper, we investigate connections between various subclasses
of harmonic univalent functions by using a convolution operator involving
the Pascal distribution series. 
\end{abstract}

\maketitle

\section{Introduction}

Let $\mathcal{H}$ denote the family of continuous complex valued harmonic
functions of the form $f=h+\overline{g}$ defined in the open unit disk $%
\mathfrak{U=}\left\{ z:\left\vert z\right\vert <1\right\} $ , where 
\begin{equation}
h(z)=z+\sum_{n=2}^{\infty }a_{n}z^{n}\text{ and }g(z)=\sum_{n=1}^{\infty
}b_{n}z^{n}\text{ \ \ }  \label{2}
\end{equation}%
are analytic in $\mathfrak{U}.$\newline
A necessary and sufficient condition for $f$ \ to be locally univalent and
sense-preserving in \emph{$\mathfrak{U}$} is that $\left\vert h^{\prime
}(z)\right\vert >\left\vert g^{\prime }(z)\right\vert $ in \emph{$\mathfrak{U}$} (see \cite{clunie}).\newline
Denote by $\mathcal{SH}$ the subclass of $\mathcal{H}$ consisting of
functions $f=h+\overline{g}$ which are harmonic, univalent and
sense-preserving in \emph{$\mathfrak{U}$} and normalized by $%
f(0)=f_{z}\left( 0\right) -1=0$. One can easily show that the
sense-preserving property implies that $\left\vert b_{1}\right\vert <1.$ The
subclass $\mathcal{SH}^{0}$ of $\mathcal{SH}$ consist of all functions in $%
\mathcal{SH}$ which have the additional property $b_{1}=0.$ Note that $%
\mathcal{SH}$ reduces to the class $\mathcal{S}$ of normalized analytic
univalent functions in \emph{$\mathfrak{U}$}, if the co-analytic part of $f$ is identically zero.\newline
A function $f\in \mathcal{SH}$ is said to be harmonic starlike of order $%
\alpha $ $\left( 0\leq \alpha <1\right) $ in \emph{$\mathfrak{U}$} if and
only if%
\begin{equation}
\Re\left\{ \frac{zf_{z}\left( z\right) -\bar{z}f_{\bar{z}}\left(
z\right) }{f\left( z\right) }\right\} >\alpha ,\text{ \ \ }(z\in \mathfrak{U)%
}  \label{3}
\end{equation}%
and is said to be harmonic convex of order $\alpha $ $\left( 0\leq \alpha
<1\right) $ in \emph{$\mathfrak{U}$} if and only if%
\begin{equation}
\Re\left\{ \frac{z^{2}f_{zz}\left( z\right) +zf_{z}\left( z\right) +%
\bar{z}^{2}f_{\bar{z}\bar{z}}\left( z\right) +\bar{z}f_{\bar{z}}\left(
z\right) }{zf_{z}\left( z\right) -\bar{z}f_{\bar{z}}\left( z\right) }%
\right\} >\alpha ,\text{ \ \ }(z\in \mathfrak{U)}.  \label{4}
\end{equation}%
These classes represented by $\mathcal{SH}^{\ast }\left( \alpha \right) $
and $\mathcal{KH}\left( \alpha \right) $, respectively, were extensively
studied by Jahangiri \cite{jahangiri}. Denote by $\mathcal{SH}^{\ast }$
and $\mathcal{KH}$ the classes $\mathcal{SH}^{\ast }(0)$ and $\mathcal{KH}%
(0),$ respectively$.$ For definitions and properties of these classes, one
may refer to \cite{silverman},\cite{silvia} or \cite{avci}.\newline
The elementary distributions such as the Poisson, the Pascal, the
Logarithmic, the Binomial have been partially studied in the Geometric
Function Theory from a theoretical point of view (see\cite{al}, \cite{por}, \cite{po}, \cite{wa}).\newline
Let us consider a non-negative discrete random variable $\mathcal{X}$ with a
Pascal probability generating function%
\begin{equation*}
\emph{P}\left( \mathcal{X}=n\right) =\binom{n+r-1}{r-1}\mathit{p}^{n}\left(
1-\mathit{p}\right) ^{r},\text{ \ \ }n\in \left\{ 0,1,2,3,...\right\} 
\end{equation*}%
where $\mathit{p},$ $r$ are called the parameters.\newline
Now we introduce a power series whose coefficients are probabilities of the
Pascal distribution, that is 
\begin{equation}
P_{p}^{r}(z)=z+\sum_{n=2}^{\infty }\binom{n+r-2}{r-1}\mathit{p}^{n-1}\left(
1-\mathit{p}\right) ^{r}z^{n}.\text{ \ \ }\left( r\geq 1,\text{ }0\leq p\leq
1,\text{ }z\in \mathfrak{U}\text{ }\right)   \label{5}
\end{equation}%
Note that, by using ratio test we conclude that the radius of convergence of
the above power series is infinity. Now, for $r,s\geq 1$ and $0\leq p,q\leq
1,$ we introduce the operator 
\begin{equation*}
P_{p,q}^{r,s}(f)(z)=P_{p}^{r}\left( z\right) \ast h\left( z\right) +%
\overline{P_{q}^{s}(z)\ast g\left( z\right) }=H\left( z\right) +\overline{%
G\left( z\right) }
\end{equation*}%
where%
\begin{eqnarray}
H(z) &=&z+\sum_{n=2}^{\infty }\binom{n+r-2}{r-1}\mathit{p}^{n-1}\left( 1-%
\mathit{p}\right) ^{r}a_{n}z^{n}  \label{6} \\
&&  \notag \\
G(z) &=&b_{1}z+\sum_{n=2}^{\infty }\binom{n+s-2}{s-1}\mathit{q}^{n-1}\left(
1-\mathit{q}\right) ^{s}b_{n}z^{n}  \notag
\end{eqnarray}%
and "$\ast $" denotes the convolution (or Hadamard product) of power series.

\section{Preliminary Lemmas}

To prove our theorems we will use the following lemmas.

\begin{lem}
(See \cite{duren}) If $f=h+\overline{g}\in \mathcal{KH}^{0}$ where $h$ and $g$
are given by $\left( \ref{2}\right) $ with $b_{1}=0,$ then%
\begin{equation*}
\left\vert a_{n}\right\vert \leq \frac{n+1}{2},\text{ \ }\left\vert
b_{n}\right\vert \leq \frac{n-1}{2}.
\end{equation*}
\end{lem}

\begin{lem}
(See \cite{jahangiri}) Let $f=h+\overline{g}$ be given by $\left( \ref{2}\right)
.$ If for some $\alpha $ $\left( 0\leq \alpha <1\right) $ and the inequality%
\begin{equation}
\sum_{n=2}^{\infty }\left( n-\alpha \right) \left\vert a_{n}\right\vert
+\sum_{n=1}^{\infty }\left( n+\alpha \right) \left\vert b_{n}\right\vert
\leq 1-\alpha   \label{7}
\end{equation}%
is hold, then $f$ is harmonic, sense-preserving, univalent in \emph{$%
\mathfrak{U}$} and $f\in \mathcal{SH}^{\ast }\left( \alpha \right) .$
\end{lem}

Define $\mathcal{TSH}^{\ast }\left( \alpha \right) =\mathcal{SH}^{\ast
}\left( \alpha \right) \cap \mathcal{T}^{2}$ and $\mathcal{TKH}\left( \alpha
\right) =\mathcal{KH}\left( \alpha \right) \cap \mathcal{T}^{1}$ where $\mathcal{%
T}^{k},~(k=1,~2)$ consisting of the functions $f=h+\overline{g}$ in SH so that $h\left(
z\right) $ and $g\left( z\right) $ are of the form%
\begin{equation}
h(z)=z-\sum_{n=2}^{\infty }\left\vert a_{n}\right\vert
z^{n},~g(z)=(-1)^{k}\sum_{n=1}^{\infty }\left\vert b_{n}\right\vert
z^{n},~|b_{1}|<1~(k=1,~2).  \label{8}
\end{equation}

\begin{rem}
(See \cite{jahangiri}) Let $f=h+\overline{g}$ be given by $\left( \ref{8}\right)
.$ Then $f\in $ $\mathcal{TSH}^{\ast }\left( \alpha \right) $ if and only if
the coefficient condition $\left( \ref{7}\right) $ is satisfied. Also, if $%
f\in \mathcal{TSH}^{\ast }\left( \alpha \right) $, then%
\begin{equation}
\left\vert a_{n}\right\vert \leq \frac{1-\alpha }{n-\alpha },\text{ \ \ }%
n\geq 2,\text{ \ \ }\left\vert b_{n}\right\vert \leq \frac{1-\alpha }{%
n+\alpha },\text{ \ \ }n\geq 1.  \label{9}
\end{equation}
\end{rem}

\begin{lem}
(See \cite{jahangiri}) Let $f=h+\overline{g}$ be given by $\left( \ref{2}\right)
.$ If for some $\alpha $ $\left( 0\leq \alpha <1\right) $ and the inequality%
\begin{equation}
\sum_{n=2}^{\infty }n\left( n-\alpha \right) \left\vert a_{n}\right\vert
+\sum_{n=1}^{\infty }n\left( n+\alpha \right) \left\vert b_{n}\right\vert
\leq 1-\alpha   \label{10}
\end{equation}%
is hold, then $f$ is harmonic, sense-preserving, univalent in \emph{$%
\mathfrak{U}$} and $f\in \mathcal{KH}\left( \alpha \right) .$
\end{lem}

\begin{rem}
(See \cite{jahangiri}) Let $f=h+\overline{g}$ be given by $\left( \ref{8}\right)
.$ Then $f\in \mathcal{TKH}\left( \alpha \right) $ if and only if the
coefficient condition $\left( \ref{10}\right) $ holds. Also, if $f\in 
\mathcal{TKH}\left( \alpha \right) $, then%
\begin{equation}
\left\vert a_{n}\right\vert \leq \frac{1-\alpha }{n\left( n-\alpha \right) },%
\text{ \ \ }n\geq 2,\text{ \ \ }\left\vert b_{n}\right\vert \leq \frac{%
1-\alpha }{n\left( n+\alpha \right) },\text{ \ \ }n\geq 1.  \label{11}
\end{equation}
\end{rem}

\begin{lem}
(See \cite{duren}) If $f=h+\overline{g}\in \mathcal{SH}^{\ast ,0}$ where $h$ and 
$g$ are given by $\left( \ref{2}\right) $ with $b_{1}=0,$ then%
\begin{equation*}
\left\vert a_{n}\right\vert \leq \frac{\left( 2n+1\right) \left( n+1\right) 
}{6},\text{ \ }\left\vert b_{n}\right\vert \leq \frac{\left( 2n-1\right)
\left( n-1\right) }{6},\text{ }n\geq 2.
\end{equation*}
\end{lem}

\section{Main Results}

\begin{thm}
Let $r,s\geq 1$ and $0\leq p,q<1.$ Also, let $f=h+\overline{g}\in \mathcal{H}
$ is given by $\left( \ref{2}\right) .$ If the inequalities  
\begin{equation}
\sum_{n=2}^{\infty }\left\vert a_{n}\right\vert +\sum_{n=1}^{\infty
}\left\vert b_{n}\right\vert \leq 1,\ \ (\left\vert b_{1}\right\vert <1)
\label{eq(i)}
\end{equation}%
and%
\begin{equation}
\left( 1-p\right) ^{r}+\left( 1-q\right) ^{s}\geq 1+\left\vert
b_{1}\right\vert +\frac{rp}{1-p}+\frac{sq}{1-q}  \label{eqii}
\end{equation}%
are hold, then $P_{p,q}^{r,s}(f)\in \mathcal{SH}^{\ast }.$
\end{thm}

\begin{proof}
Note that $P_{p,q}^{r,s}(f)=H\left( z\right) +\overline{G\left( z\right) },$
where $H\left( z\right) $ and $G\left( z\right) $ are given by $\left( \ref%
{6}\right) .$ To prove that  $P_{p,q}^{r,s}(f)$ is locally univalent and
sense-preserving it suffices to prove that \ $\left\vert H^{\prime
}(z)\right\vert -\left\vert G^{\prime }(z)\right\vert >0$ in \emph{$%
\mathfrak{U}$}$.$ Using (\ref{eq(i)}), we compute
\begin{scriptsize}
\begin{eqnarray*}
\left\vert H^{\prime }(z)\right\vert -\left\vert G^{\prime }(z)\right\vert 
&>&1-\sum_{n=2}^{\infty }n\binom{n+r-2}{r-1}\mathit{p}^{n-1}\left( 1-\mathit{%
p}\right) ^{r} \\
&&-\left\vert b_{1}\right\vert -\sum_{n=2}^{\infty }n\binom{n+s-2}{s-1}%
\mathit{q}^{n-1}\left( 1-\mathit{q}\right) ^{s} \\
&=&1-\left\vert b_{1}\right\vert -\sum_{n=2}^{\infty }\left( n-1+1\right) 
\binom{n+r-2}{r-1}\mathit{p}^{n-1}\left( 1-\mathit{p}\right) ^{r} \\
&&-\sum_{n=2}^{\infty }\left( n-1+1\right) \binom{n+s-2}{s-1}\mathit{q}%
^{n-1}\left( 1-\mathit{q}\right) ^{s} \\
&=&1-\left\vert b_{1}\right\vert -rp\left( 1-\mathit{p}\right)
^{r}\sum_{n=2}^{\infty }\binom{n+r-2}{r}\mathit{p}^{n-2} \\
&&-\left( 1-\mathit{p}\right) ^{r}\sum_{n=2}^{\infty }\binom{n+r-2}{r-1}%
\mathit{p}^{n-1}-sq\left( 1-\mathit{q}\right) ^{s}\sum_{n=2}^{\infty }\binom{%
n+s-2}{s}\mathit{q}^{n-2} \\
&&-\left( 1-\mathit{q}\right) ^{s}\sum_{n=2}^{\infty }\binom{n+s-2}{s-1}%
\mathit{q}^{n-1}
\end{eqnarray*}
\end{scriptsize}
\begin{eqnarray*}
&=&1-\left\vert b_{1}\right\vert -rp\left( 1-\mathit{p}\right)
^{r}\sum_{n=0}^{\infty }\binom{n+r}{r}\mathit{p}^{n} \\
&&-\left( 1-\mathit{p}\right) ^{r}\sum_{n=0}^{\infty }\binom{n+r-1}{r-1}%
\mathit{p}^{n}+\left( 1-\mathit{p}\right) ^{r} \\
&&-sq\left( 1-\mathit{q}\right) ^{s}\sum_{n=0}^{\infty }\binom{n+s}{s}%
\mathit{q}^{n} \\
&&-\left( 1-\mathit{q}\right) ^{s}\sum_{n=0}^{\infty }\binom{n+s-1}{s-1}%
\mathit{q}^{n}+\left( 1-\mathit{q}\right) ^{s} \\
&=&\left( 1-\mathit{p}\right) ^{r}+\left( 1-\mathit{q}\right)
^{s}-1-\left\vert b_{1}\right\vert -\frac{rp}{1-p}-\frac{sq}{1-q}\geq 0.
\end{eqnarray*}%
To prove $P_{p,q}^{r,s}(f)$ is univalent in \emph{$\mathfrak{U}$}, referring
Theorem 1 in \cite{jahangiri}, for $z_{1}\neq z_{2}$ in \emph{$\mathfrak{U}
$}$,$ we need to show that 
\begin{equation}
\Re\frac{P_{p,q}^{r,s}(f)\left( z_{2}\right) -P_{p,q}^{r,s}(f)\left(
z_{1}\right) }{z_{2}-z_{1}}>\int_{0}^{1}\left( \Re(H^{\prime }\left(
z\left( t\right) \right)) -\left\vert G^{\prime }\left( z\left( t\right)
\right) \right\vert \right) dt.  \label{12}
\end{equation}%
By (\ref{eq(i)}), we have 
\begin{eqnarray*}
\Re(H^{\prime }\left( z\left( t\right) \right)) -\left\vert G^{\prime
}\left( z\left( t\right) \right) \right\vert  &>&1-\sum_{n=2}^{\infty }n%
\binom{n+r-2}{r-1}\mathit{p}^{n-1}\left( 1-\mathit{p}\right) ^{r} \\
&&-\left\vert b_{1}\right\vert -\sum_{n=2}^{\infty }n\binom{n+s-2}{s-1}%
\mathit{q}^{n-1}\left( 1-\mathit{q}\right) ^{s}.
\end{eqnarray*}%
Using\ (\ref{eqii}), we obtain that the inequality above is nonnegative.
Therefore, from the inequality $\left( \ref{12}\right) $ we have 
\begin{equation*}
\Re\frac{P_{p,q}^{r,s}(f)\left( z_{2}\right) -P_{p,q}^{r,s}(f)\left(
z_{1}\right) }{z_{2}-z_{1}}>0.
\end{equation*}%
Hence univalency of $P_{p,q}^{r,s}(f)$ is proved. \newline
In order to show that $P_{p,q}^{r,s}(f)\in \mathcal{SH}^{\ast },$ we need to
prove $\Phi _{1}\leq 1,$ by Lemma 2.2, where%
\begin{equation*}
\Phi _{1}=\sum_{n=2}^{\infty }n\binom{n+r-2}{r-1}\mathit{p}^{n-1}\left( 1-%
\mathit{p}\right) ^{r}\left\vert a_{n}\right\vert +\left\vert
b_{1}\right\vert +\sum_{n=2}^{\infty }n\binom{n+s-2}{s-1}\mathit{q}%
^{n-1}\left( 1-\mathit{q}\right) ^{s}\left\vert b_{n}\right\vert .
\end{equation*}%
Since $\left\vert a_{n}\right\vert \leq 1,$ $\left\vert b_{n}\right\vert
\leq 1,$ $\forall n\geq 2$ because of (\ref{eq(i)}), we have%
\begin{eqnarray*}
\Phi _{1} &\leq &rp\left( 1-\mathit{p}\right) ^{r}\sum_{n=0}^{\infty }\binom{%
n+r}{r}\mathit{p}^{n}+\left( 1-\mathit{p}\right) ^{r}\sum_{n=0}^{\infty }%
\binom{n+r-1}{r-1}\mathit{p}^{n} \\
&&-\left( 1-\mathit{p}\right) ^{r}+\left\vert b_{1}\right\vert +sq\left( 1-%
\mathit{q}\right) ^{s}\sum_{n=0}^{\infty }\binom{n+s}{s}\mathit{q}^{n} \\
&&+\left( 1-\mathit{q}\right) ^{s}\sum_{n=0}^{\infty }\binom{n+s-1}{s-1}%
\mathit{q}^{n}-\left( 1-\mathit{q}\right) ^{s} \\
&=&\left\vert b_{1}\right\vert +\frac{rp}{1-p}+1-\left( 1-\mathit{p}\right)
^{r}+\frac{sq}{1-q}+1-\left( 1-\mathit{q}\right) ^{s} \\
&\leq &1
\end{eqnarray*}%
from (\ref{eqii}). Thus proof of Theorem 3.1 is completed.
\end{proof}

\begin{thm}
Let $0\leq \alpha <1$, $r,s\geq 1$ and $0\leq p,q<1.$ If the inequality%
\begin{eqnarray*}
&&\frac{r\left( r+1\right) \mathit{p}^{2}}{\left( 1-\mathit{p}\right) ^{2}}+%
\frac{\left( 4-\alpha \right) r\mathit{p}}{1-\mathit{p}}+\frac{s\left(
s+1\right) \mathit{q}^{2}}{\left( 1-\mathit{q}\right) ^{2}}+\frac{\left(
2+\alpha \right) s\mathit{q}}{1-\mathit{q}} \\
&&\leq 2\left( 1-\alpha \right) (1-p)^{r}
\end{eqnarray*}%
is hold, then $P_{p,q}^{r,s}\left( \mathcal{KH}^{0}\right) \subset \mathcal{SH}^{\ast ,0}\left( \alpha \right)$.
\end{thm}

\begin{proof}
Suppose that $f=h+\overline{g}\in \mathcal{KH}^{0}$ where $h$ and $g$ are
given by $\left( \ref{2}\right) $ with $b_{1}=0.$ It suffices to show that $%
P_{p,q}^{r,s}(f)=H+\overline{G}\in \mathcal{SH}^{\ast ,0}\left( \alpha
\right) $ where $H$ and $G$ are given by $\left( \ref{6}\right) $ with $%
b_{1}=0$ in \emph{$\mathfrak{U}$}$.$ Using Lemma 2.2, we need to prove that $%
\Phi _{2}\leq 1-\alpha ,$ where%
\begin{eqnarray}
\Phi _{2} &=&\sum_{n=2}^{\infty }\left( n-\alpha \right) \binom{n+r-2}{r-1}%
\left( 1-\mathit{p}\right) ^{r}\mathit{p}^{n-1}\left\vert a_{n}\right\vert 
\label{13} \\
&&+\sum_{n=2}^{\infty }\left( n+\alpha \right) \binom{n+s-2}{s-1}\left( 1-%
\mathit{q}\right) ^{s}\mathit{q}^{n-1}\left\vert b_{n}\right\vert .
\end{eqnarray}%
Using Lemma 2.1, we compute
\begin{scriptsize}
\begin{eqnarray*}
\Phi _{2} &\leq &\frac{1}{2}\left\{ \sum_{n=2}^{\infty }\left( n-\alpha
\right) \left( n+1\right) \binom{n+r-2}{r-1}\left( 1-\mathit{p}\right) ^{r}%
\mathit{p}^{n-1}\right.  \\
&&\left. +\sum_{n=2}^{\infty }\left( n+\alpha \right) \left( n-1\right) 
\binom{n+s-2}{s-1}\left( 1-\mathit{q}\right) ^{s}\mathit{q}^{n-1}\right\}  \\
&=&\frac{1}{2}\left\{ \sum_{n=2}^{\infty }\left[ \left( n-1\right) \left(
n-2\right) +\left( 4-\alpha \right) \left( n-1\right) +2\left( 1-\alpha
\right) \right] \binom{n+r-2}{r-1}\left( 1-\mathit{p}\right) ^{r}\mathit{p}%
^{n-1}\right.  \\
&&\left. +\sum_{n=2}^{\infty }\left[ \left( n-1\right) \left( n-2\right)
+\left( 2+\alpha \right) \left( n-1\right) \right] \binom{n+s-2}{s-1}\left(
1-\mathit{q}\right) ^{s}\mathit{q}^{n-1}\right\}  \\
&=&\frac{1}{2}\left\{ r\left( r+1\right) \mathit{p}^{2}\left( 1-\mathit{p}%
\right) ^{r}\sum_{n=3}^{\infty }\binom{n+r-2}{r+1}\mathit{p}^{n-3}\right.  \\
&&\left. +\left( 4-\alpha \right) r\mathit{p}\left( 1-\mathit{p}\right)
^{r}\sum_{n=2}^{\infty }\binom{n+r-2}{r}\mathit{p}^{n-2}\right.  \\
&&\left. +2\left( 1-\alpha \right) \left( 1-\mathit{p}\right)
^{r}\sum_{n=2}^{\infty }\binom{n+r-2}{r-1}\mathit{p}^{n-2}\right.  \\
&&\left. +s\left( s+1\right) \mathit{q}^{2}\left( 1-\mathit{q}\right)
^{s}\sum_{n=3}^{\infty }\binom{n+s-2}{s+1}\mathit{q}^{n-3}\right.  \\
&&\left. +\left( 2+\alpha \right) s\mathit{q}\left( 1-\mathit{q}\right)
^{s}\sum_{n=2}^{\infty }\binom{n+s-2}{s}\mathit{q}^{n-2}\right\} 
\end{eqnarray*}
\end{scriptsize}

\begin{eqnarray*}
&=&\frac{1}{2}\left\{ r\left( r+1\right) \mathit{p}^{2}\left( 1-\mathit{p}%
\right) ^{r}\sum_{n=0}^{\infty }\binom{n+r+1}{r+1}\mathit{p}^{n}\right.  \\
&&\left. +\left( 4-\alpha \right) r\mathit{p}\left( 1-\mathit{p}\right)
^{r}\sum_{n=0}^{\infty }\binom{n+r}{r}\mathit{p}^{n}\right.  \\
&&\left. +2\left( 1-\alpha \right) \left( 1-\mathit{p}\right)
^{r}\sum_{n=0}^{\infty }\binom{n+r-1}{r-1}\mathit{p}^{n}-2\left( 1-\alpha
\right) \left( 1-\mathit{p}\right) ^{r}\right.  \\
&&\left. +s\left( s+1\right) \mathit{q}^{2}\left( 1-\mathit{q}\right)
^{s}\sum_{n=0}^{\infty }\binom{n+s+1}{s+1}\mathit{q}^{n}\right.  \\
&&\left. +\left( 2+\alpha \right) s\mathit{q}\left( 1-\mathit{q}\right)
^{s}\sum_{n=0}^{\infty }\binom{n+s}{s}\mathit{q}^{n}\right\}  \\
&=&\frac{1}{2}\left\{ \frac{r\left( r+1\right) \mathit{p}^{2}}{\left( 1-%
\mathit{p}\right) ^{2}}+\frac{\left( 4-\alpha \right) r\mathit{p}}{1-\mathit{%
p}}+2\left( 1-\alpha \right) \right.  \\
&&\left. -2\left( 1-\alpha \right) \left( 1-\mathit{p}\right) ^{r}+\frac{%
s\left( s+1\right) \mathit{q}^{2}}{\left( 1-\mathit{q}\right) ^{2}}+\frac{%
\left( 2+\alpha \right) s\mathit{q}}{1-\mathit{q}}\right\} .
\end{eqnarray*}%
The last expression is bounded above by $\left( 1-\alpha \right) $ by the
given condition. Thus the proof of Theorem 3.2 is completed.
\end{proof}

\begin{thm}
Suppose $0\leq \alpha <1$, $r,s\geq 1$ and $0\leq p,q<1.$ If the inequality%
\begin{eqnarray}
&&\frac{2r\left( r+1\right) \left( r+2\right) \mathit{p}^{3}}{\left( 1-%
\mathit{p}\right) ^{3}}+\frac{\left( 15-2\alpha \right) r\left( r+1\right) 
\mathit{p}^{2}}{\left( 1-\mathit{p}\right) ^{2}}+\frac{\left( 24-9\alpha
\right) r\mathit{p}}{1-\mathit{p}}  \label{eqtheo9} \\
&&+\frac{2s\left( s+1\right) \left( s+2\right) \mathit{q}^{3}}{\left( 1-%
\mathit{q}\right) ^{3}}+\frac{\left( 9+2\alpha \right) s\left( s+1\right) 
\mathit{q}^{2}}{\left( 1-\mathit{q}\right) ^{2}}+\frac{\left( 6+3\alpha
\right) s\mathit{q}}{1-\mathit{q}}  \notag \\
&\leq &6\left( 1-\alpha \right) \left( 1-\mathit{p}\right) ^{r}  \notag
\end{eqnarray}%
is hold then $P_{p,q}^{r,s}\left( \mathcal{SH}^{\ast ,0}\left( \alpha
\right) \right) \subset \mathcal{SH}^{\ast ,0}\left( \alpha \right) .$
\end{thm}

\begin{proof}
Suppose $f=h+\overline{g}\in \mathcal{SH}^{\ast .0}\left( \alpha \right) $
where $h$ and $g$ are given by $\left( \ref{2}\right) $ with $b_{1}=0.$ It
suffices to show that $P_{p,q}^{r,s}(f)=H+\overline{G}\in \mathcal{SH}^{\ast
,0}\left( \alpha \right) $ where $H$ and $G$ are given by $\left( \ref{6}%
\right) $ with $b_{1}=0.$ By Lemma 2.2, we need to prove that $\Phi _{2}\leq
1-\alpha ,$ where%
\begin{eqnarray*}
\Phi _{2} &=&\sum_{n=2}^{\infty }\left( n-\alpha \right) \binom{n+r-2}{r-1}%
\left( 1-\mathit{p}\right) ^{r}\mathit{p}^{n-1}\left\vert a_{n}\right\vert 
\\
&&+\sum_{n=2}^{\infty }\left( n+\alpha \right) \binom{n+s-2}{s-1}\left( 1-%
\mathit{q}\right) ^{s}\mathit{q}^{n-1}\left\vert b_{n}\right\vert .
\end{eqnarray*}%
Using Lemma 2.6, we have 
\begin{eqnarray*}
\Phi _{2} &\leq &\frac{1}{6}\left\{ \sum_{n=2}^{\infty }\left( n-\alpha
\right) \left( 2n+1\right) \left( n+1\right) \binom{n+r-2}{r-1}\left( 1-%
\mathit{p}\right) ^{r}\mathit{p}^{n-1}\right.  \\
&&\left. +\sum_{n=2}^{\infty }\left( n+\alpha \right) \left( 2n-1\right)
\left( n-1\right) \binom{n+s-2}{s-1}\left( 1-\mathit{q}\right) ^{s}\mathit{q}%
^{n-1}\right\}  \\
&& \\
&=&\frac{1}{6}\left\{ 2\sum_{n=2}^{\infty }\binom{n+r-2}{r-1}\left(
n-1\right) \left( n-2\right) \left( n-3\right) \left( 1-\mathit{p}\right)
^{r}\mathit{p}^{n-1}\right.  \\
&&+\left( 15-2\alpha \right) \sum_{n=2}^{\infty }\binom{n+r-2}{r-1}\left(
n-1\right) \left( n-2\right) \left( 1-\mathit{p}\right) ^{r}\mathit{p}^{n-1}
\\
&&+\left( 24-9\alpha \right) \sum_{n=2}^{\infty }\binom{n+r-2}{r-1}\left(
n-1\right) \left( 1-\mathit{p}\right) ^{r}\mathit{p}^{n-1} \\
&&+6\left( 1-\alpha \right) \sum_{n=2}^{\infty }\binom{n+r-2}{r-1}\left( 1-%
\mathit{p}\right) ^{r}\mathit{p}^{n-1} \\
&&+2\sum_{n=2}^{\infty }\binom{n+s-2}{s-1}\left( n-1\right) \left(
n-2\right) \left( n-3\right) \left( 1-\mathit{q}\right) ^{s}\mathit{q}^{n-1}
\\
&&+\left( 9+2\alpha \right) \sum_{n=2}^{\infty }\binom{n+s-2}{s-1}\left(
n-1\right) \left( n-2\right) \left( 1-\mathit{q}\right) ^{s}\mathit{q}^{n-1}
\\
&&\left. +\left( 6+3\alpha \right) \sum_{n=2}^{\infty }\binom{n+s-2}{s-1}%
\left( n-1\right) \left( 1-\mathit{q}\right) ^{s}\mathit{q}^{n-1}\right\}  \\
&=&\frac{1}{6}\left\{ 2r\left( r+1\right) \left( r+2\right) \mathit{p}%
^{3}\left( 1-\mathit{p}\right) ^{r}\sum_{n=4}^{\infty }\binom{n+r-2}{r+2}%
\mathit{p}^{n-4}\right.  \\
&&+\left( 15-2\alpha \right) r\left( r+1\right) \mathit{p}^{2}\left( 1-%
\mathit{p}\right) ^{r}\sum_{n=3}^{\infty }\binom{n+r-2}{r+1}\mathit{p}^{n-3}
\\
&&+\left( 24-9\alpha \right) r\mathit{p}\left( 1-\mathit{p}\right)
^{r}\sum_{n=2}^{\infty }\binom{n+r-2}{r}\mathit{p}^{n-2} \\
&&+6\left( 1-\alpha \right) \sum_{n=2}^{\infty }\binom{n+r-2}{r-1}\left( 1-%
\mathit{p}\right) ^{r}\mathit{p}^{n-1} \\
&&+2s\left( s+1\right) \left( s+2\right) \mathit{q}^{3}\left( 1-\mathit{q}%
\right) ^{s}\sum_{n=4}^{\infty }\binom{n+s-2}{s+2}\mathit{q}^{n-4} \\
&&+\left( 9+2\alpha \right) s\left( s+1\right) \mathit{q}^{2}\left( 1-%
\mathit{q}\right) ^{s}\sum_{n=3}^{\infty }\binom{n+s-2}{s+1}\mathit{q}^{n-3}
\\
&&\left. +\left( 6+3\alpha \right) s\mathit{q}\left( 1-\mathit{q}\right)
^{s}\sum_{n=2}^{\infty }\binom{n+s-2}{s}\mathit{q}^{n-2}\right\} 
\end{eqnarray*}%
\begin{eqnarray*}
&=&\frac{1}{6}\left\{ 2r\left( r+1\right) \left( r+2\right) \mathit{p}%
^{3}\left( 1-\mathit{p}\right) ^{r}\sum_{n=0}^{\infty }\binom{n+r+2}{r+2}%
\mathit{p}^{n}\right.  \\
&&+\left( 15-2\alpha \right) r\left( r+1\right) \mathit{p}^{2}\left( 1-%
\mathit{p}\right) ^{r}\sum_{n=0}^{\infty }\binom{n+r+1}{r+1}\mathit{p}^{n} \\
&&+\left( 24-9\alpha \right) r\mathit{p}\left( 1-\mathit{p}\right)
^{r}\sum_{n=0}^{\infty }\binom{n+r}{r}\mathit{p}^{n} \\
&&+6\left( 1-\alpha \right) \left( 1-\mathit{p}\right)
^{r}\sum_{n=0}^{\infty }\binom{n+r-1}{r-1}\mathit{p}^{n}-6\left( 1-\alpha
\right) \left( 1-\mathit{p}\right) ^{r} \\
&&+2s\left( s+1\right) \left( s+2\right) \mathit{q}^{3}\left( 1-\mathit{q}%
\right) ^{s}\sum_{n=0}^{\infty }\binom{n+s+2}{s+2}\mathit{q}^{n} \\
&&+\left( 9+2\alpha \right) s\left( s+1\right) \mathit{q}^{2}\left( 1-%
\mathit{q}\right) ^{s}\sum_{n=0}^{\infty }\binom{n+s+1}{s+1}\mathit{q}^{n} \\
&&\left. +\left( 6+3\alpha \right) s\mathit{q}\left( 1-\mathit{q}\right)
^{s}\sum_{n=0}^{\infty }\binom{n+s}{s}\mathit{q}^{n}\right\} 
\end{eqnarray*}

\begin{eqnarray*}
&=&\frac{1}{6}\left \{ \frac{2r\left( r+1\right) \left( r+2\right) \mathit{p}%
^{3}}{\left( 1-\mathit{p}\right) ^{3}}+\frac{\left( 15-2\alpha \right)
r\left( r+1\right) \mathit{p}^{2}}{\left( 1-\mathit{p}\right) ^{2}}\right. \\
&&+\frac{\left( 24-9\alpha \right) r\mathit{p}}{1-\mathit{p}}+6\left(
1-\alpha \right) -6\left( 1-\alpha \right) \left( 1-\mathit{p}\right) ^{r} \\
&&\left. +\frac{2s\left( s+1\right) \left( s+2\right) \mathit{q}^{3}}{\left(
1-\mathit{q}\right) ^{3}}+\frac{\left( 9+2\alpha \right) s\left( s+1\right) 
\mathit{q}^{2}}{\left( 1-\mathit{q}\right) ^{2}}+\frac{\left( 6+3\alpha
\right) s\mathit{q}}{1-\mathit{q}}\right \} \\
&\leq &1-\alpha.
\end{eqnarray*}%
\end{proof}

\begin{thm}
Let$0\leq \alpha <1$, $r,s\geq 1$ and $0\leq p,q<1.$ If the inequality%
\begin{equation*}
\left( 1-\mathit{p}\right) ^{r}+\left( 1-\mathit{q}\right) ^{s}\geq 1+\frac{%
\left( 1+\alpha \right) \left\vert b_{1}\right\vert }{\left( 1-\alpha
\right) }
\end{equation*}%
is hold, then $P_{p,q}^{r,s}\left( \mathcal{TSH}^{\ast }\left( \alpha
\right) \right) \subset \mathcal{TSH}^{\ast }\left( \alpha \right) .$
\end{thm}

\begin{proof}
Suppose $f=h+\overline{g}\in T\mathcal{SH}^{\ast }\left( \alpha \right) $
where $h$ and $g$ are given by $\left( \ref{8}\right) $ with $b_{1}=0.$ We
need to prove that the operator 
\begin{eqnarray*}
P_{p,q}^{r,s}(f)\left( z\right)  &=&z-\sum_{n=2}^{\infty }\binom{n+r-2}{r-1}%
\left( 1-\mathit{p}\right) ^{r}\mathit{p}^{n-1}\left\vert a_{n}\right\vert
z^{n} \\
&&\left\vert b_{1}\right\vert \overline{z}+\sum_{n=2}^{\infty }\binom{n+s-2}{%
s-1}\left( 1-\mathit{q}\right) ^{s}\mathit{q}^{n-1}\left\vert
b_{n}\right\vert \overline{z}^{n}
\end{eqnarray*}%
is in $T\mathcal{SH}^{\ast }\left( \alpha \right) $ if and only if $\Phi
_{3}\leq 1-\alpha ,$ where%
\begin{eqnarray*}
\Phi _{3} &=&\sum_{n=2}^{\infty }\left( n-\alpha \right) \binom{n+r-2}{r-1}%
\left( 1-\mathit{p}\right) ^{r}\mathit{p}^{n-1}\left\vert a_{n}\right\vert 
\\
&&+\left( 1+\alpha \right) \left\vert b_{1}\right\vert +\sum_{n=2}^{\infty
}\left( n+\alpha \right) \binom{n+s-2}{s-1}\left( 1-\mathit{q}\right) ^{s}%
\mathit{q}^{n-1}\left\vert b_{n}\right\vert .
\end{eqnarray*}%
By Remark 2.3, we obtain
\begin{eqnarray*}
\Phi _{3} &\leq &\left( 1-\alpha \right) \left\{ \sum_{n=2}^{\infty }\binom{%
n+r-2}{r-1}\left( 1-\mathit{p}\right) ^{r}\mathit{p}^{n-1}\right.  \\
&&\left. +\sum_{n=1}^{\infty }\binom{n+s-2}{s-1}\left( 1-\mathit{q}\right)
^{s}\mathit{q}^{n-1}\right\} +\left( 1+\alpha \right) \left\vert
b_{1}\right\vert  \\
&=&\left( 1-\alpha \right) \left\{ \left( 1-\mathit{p}\right)
^{r}\sum_{n=0}^{\infty }\binom{n+r-1}{r-1}\mathit{p}^{n}-\left( 1-\mathit{p}%
\right) ^{r}\right.  \\
&&\left. +\left( 1-\mathit{q}\right) ^{s}\sum_{n=0}^{\infty }\binom{n+s-1}{%
s-1}\mathit{q}^{n}-\left( 1-\mathit{q}\right) ^{s}\right\} +\left( 1+\alpha
\right) \left\vert b_{1}\right\vert  \\
&=&\left( 1-\alpha \right) \left\{ 2-\left( 1-\mathit{p}\right) ^{r}-\left(
1-\mathit{q}\right) ^{s}\right\} +\left( 1+\alpha \right) \left\vert
b_{1}\right\vert  \\
&\leq &1-\alpha.
\end{eqnarray*}%
Thus the proof of the theorem is completed.
\end{proof}

We next explore a sufficient condition which guarantees that $P_{p,q}^{r,s}$
maps $\mathcal{KH}^{0}$ into $\mathcal{KH}\left( \alpha \right) .$

\begin{thm}
Suppose $0\leq \alpha <1$, $r,s\geq 1$ and $0\leq p,q<1.$ If the inequality%
\begin{eqnarray*}
&&\frac{r\left( r+1\right) \left( r+2\right) \mathit{p}^{3}}{\left( 1-%
\mathit{p}\right) ^{3}}+\frac{\left( 7-\alpha \right) r\left( r+1\right) 
\mathit{p}^{2}}{\left( 1-\mathit{p}\right) ^{2}}+\frac{\left( 10-4\alpha
\right) r\mathit{p}}{1-\mathit{p}} \\
&&+\frac{s\left( s+1\right) \left( s+2\right) \mathit{q}^{3}}{\left( 1-%
\mathit{q}\right) ^{3}}+\frac{\left( 5+\alpha \right) s\left( s+1\right) 
\mathit{q}^{2}}{\left( 1-\mathit{q}\right) ^{2}}+\frac{\left( 4+2\alpha
\right) s\mathit{q}}{1-\mathit{q}} \\
&\leq &2\left( 1-\alpha \right) \left( 1-\mathit{p}\right) ^{r}
\end{eqnarray*}%
is hold, then $P_{p,q}^{r,s}\left( \mathcal{KH}^{0}\right) \subset \mathcal{%
KH}^{0}\left( \alpha \right) .$
\end{thm}

\begin{proof}
Let $f=h+\overline{g}\in \mathcal{KH}^{0}$ where $h$ and $g$ are given by $%
\left( \ref{2}\right) $ with $b_{1}=0.$ It suffices to show that $%
P_{p,q}^{r,s}(f)=H+\overline{G}\in \mathcal{KH}^{0}\left( \alpha \right) $
where $H$ and $G$ are given by $\left( \ref{6}\right) $ with $b_{1}=0.$
Referring Lemma 2.1, we need to prove that $\Phi _{4}\leq 1-\alpha ,$ where%
\begin{eqnarray*}
\Phi _{4} &=&\sum_{n=2}^{\infty }n\left( n-\alpha \right) \binom{n+r-2}{r-1}%
\left( 1-\mathit{p}\right) ^{r}\mathit{p}^{n-1}\left\vert a_{n}\right\vert 
\\
&&+\sum_{n=2}^{\infty }n\left( n+\alpha \right) \binom{n+s-2}{s-1}\left( 1-%
\mathit{q}\right) ^{s}\mathit{q}^{n-1}\left\vert b_{n}\right\vert .
\end{eqnarray*}%
Hence, 
\begin{eqnarray*}
\Phi _{4} &\leq &\frac{1}{2}\left\{ \sum_{n=2}^{\infty }\left( n-1\right)
\left( n-2\right) \left( n-3\right) \binom{n+r-2}{r-1}\left( 1-\mathit{p}%
\right) ^{r}\mathit{p}^{n-1}\right.  \\
&&+\sum_{n=2}^{\infty }\left( 7-\alpha \right) \left( n-1\right) \left(
n-2\right) \binom{n+r-2}{r-1}\left( 1-\mathit{p}\right) ^{r}\mathit{p}^{n-1}
\\
&&+\sum_{n=2}^{\infty }\left( 10-4\alpha \right) \left( n-1\right)
+\sum_{n=2}^{\infty }2\left( 1-\alpha \right) \binom{n+r-2}{r-1}\left( 1-%
\mathit{p}\right) ^{r}\mathit{p}^{n-1} \\
&&+\sum_{n=2}^{\infty }\left( n-1\right) \left( n-2\right) \left( n-3\right) 
\binom{n+s-2}{s-1}\left( 1-\mathit{q}\right) ^{s}\mathit{q}^{n-1} \\
&&+\sum_{n=2}^{\infty }\left( 5+\alpha \right) \left( n-1\right) \left(
n-2\right) \binom{n+s-2}{s-1}\left( 1-\mathit{q}\right) ^{s}\mathit{q}^{n-1}
\\
&&\left. +\sum_{n=2}^{\infty }\left( 4+2\alpha \right) \left( n-1\right) 
\binom{n+s-2}{s-1}\left( 1-\mathit{q}\right) ^{s}\mathit{q}^{n-1}\right\}  \\
&=&\frac{1}{2}\left\{ \frac{r\left( r+1\right) \left( r+2\right) \mathit{p}%
^{3}}{\left( 1-\mathit{p}\right) ^{3}}+\frac{\left( 7-\alpha \right) r\left(
r+1\right) \mathit{p}^{2}}{\left( 1-\mathit{p}\right) ^{2}}+\frac{\left(
10-4\alpha \right) r\mathit{p}}{1-\mathit{p}}\right.  \\
&&+2\left( 1-\alpha \right) -2\left( 1-\alpha \right) \left( 1-\mathit{p}%
\right) ^{r} \\
&&\left. +\frac{s\left( s+1\right) \left( s+2\right) \mathit{q}^{3}}{\left(
1-\mathit{q}\right) ^{3}}+\frac{\left( 5+\alpha \right) s\left( s+1\right) 
\mathit{q}^{2}}{\left( 1-\mathit{q}\right) ^{2}}+\frac{\left( 4+2\alpha
\right) s\mathit{q}}{1-\mathit{q}}\right\}  \\
&\leq &1-\alpha.
\end{eqnarray*}%
\end{proof}

The proofs of following theorems are similar to previous theorems so we omit
them.

\begin{thm}
If $0\leq \alpha <1$, $r,s\geq 1$ and $0\leq p,q<1$ then $%
P_{p,q}^{r,s}\left( \mathcal{TSH}^{\ast }\left( \alpha \right) \right)
\subset \mathcal{TKH}\left( \alpha \right) $ if and only if the inequality%
\begin{equation}
\left( 1-\mathit{p}\right) ^{r}+\left( 1-\mathit{q}\right) ^{s}\geq 1+\frac{r%
\mathit{p}}{1-\mathit{p}}+\frac{s\mathit{q}}{1-\mathit{q}}+\frac{\left(
1+\alpha \right) }{\left( 1-\alpha \right) }\left\vert b_{1}\right\vert 
\label{eqtheo12}
\end{equation}%
is hold.
\end{thm}

\begin{thm}
If $0\leq \alpha <1$, $r,s\geq 1$ and $0\leq p,q<1$ then $%
P_{p,q}^{r,s}\left( \mathcal{TKH}\left( \alpha \right) \right) \subset 
\mathcal{TKH}\left( \alpha \right) $ if and only if the inequality%
\begin{equation*}
\left( 1-\mathit{p}\right) ^{r}+\left( 1-\mathit{q}\right) ^{s}\geq 1+\frac{%
\left( 1+\alpha \right) \left\vert b_{1}\right\vert }{\left( 1-\alpha
\right) }
\end{equation*}%
is hold.
\end{thm}

\begin{ex}
Consider the harmonic function $f(z)=z+\frac{1}{5}\overline{z}^{2}.$ If we
take $r=2,$ $s=2,$ $p=0.01$ and $q=0.01$ then from (\ref{6}), we have 
\begin{equation*}
P_{0.01,\text{ }0.01}^{2,2}(f)(z)=z+0.0039\overline{z}^{2}.
\end{equation*}
Then we get the following results:\newline
$(i)$ since condition (\ref{eq(i)}) is satisfied, by Theorem 3.1, $P_{0.01,%
\text{ }0.01}^{2,2}(f)\in \mathcal{SH}^{\ast },$\newline
$(ii)$ since condition (\ref{eqtheo9}) is satisfied $f\in \mathcal{SH}^{\ast
}(\frac{1}{2})$, by Theorem 3.3, $P_{0.01,\text{ }0.01}^{2,2}(f)\in \mathcal{SH%
}^{\ast }(\frac{1}{2}),$\newline
$(iii)$ since condition (\ref{eqtheo12}) is satisfied $f\in \mathcal{TSH}%
^{\ast }(\frac{1}{2}),$ by Theorem 3.6, $P_{0.01,\text{ }0.01}^{2,2}(f)\in 
\mathcal{TKH}^{0}(\frac{1}{2}).$\newline
Images of concentric circles inside $\mathfrak{U}$ under the functions $f$
and $P_{0.01,\text{ }0.01}^{2,2}(f)$ are shown in Figures 1 and 2.%

\begin{figure}[h]
	\begin{minipage}[h]{0.40\linewidth}
		\centering
		\includegraphics[width=6cm, height=6cm]{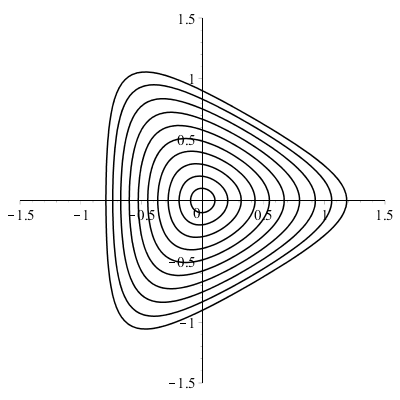}
		\caption{{\small Image of $f(\mathfrak{U})$}}
	\end{minipage}
	\hspace{0.1\linewidth}
	\begin{minipage}[h]{0.40\linewidth}
		\centering
		\includegraphics[width=6cm, height=6cm]{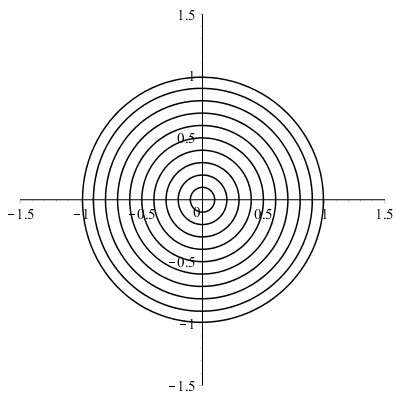}
		\caption{     {\small Image of $P_{0.01,0.01}^{2,2}(f(\mathfrak{U}))$} }
	
	\end{minipage}
\end{figure}
\end{ex}

\begin{ex}
Consider the harmonic right half plane mapping $f_{0}(z)=\frac{z-\frac{1}{2}%
z^{2}}{(1-z)^{2}}+\frac{-\frac{1}{2}\overline{z}^{2}}{(1-\overline{z})^{2}}%
\in \mathcal{K}\mathcal{H}^{0}.$ If we take $r=2,$ $s=2,$ $p=0.01$ and $%
q=0.01$ then from (\ref{6}), we have 
\begin{eqnarray*}
P_{0.01,\text{ }0.01}^{2,2}(f_{0})(z) &=&z+\sum_{n=2}^{\infty }\frac{n(n+1)}{%
2}\text{ }(0.01)^{n-1}(0.99)^{2}z^{n} \\
&&+\sum_{n=2}^{\infty }\frac{n(-n+1)}{2}(0.01)^{n-1}(0.99)^{2}%
\overline{z}^{n}.
\end{eqnarray*}%
Then, according to the Theorem 3.5, $P_{0.01,\text{ }0.01}^{2,2}(f_{0})(z)\in 
\mathcal{K}\mathcal{H}^{0}(\alpha )$ for $0\leq \alpha <1.$ Images of
concentric circles inside $\mathfrak{U}$ under the functions $f_{0}$ and $%
P_{0.01,\text{ }0.01}^{2,2}(f_{0})$ are shown in Figures 3 and 4.
\begin{figure}[h]
	\begin{minipage}[h]{0.40\linewidth}
		\centering
		\includegraphics[width=6cm, height=6cm]{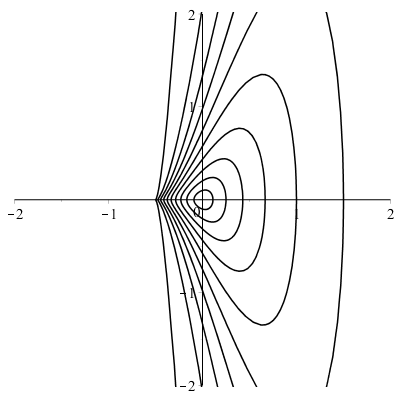}
		\caption{{\small Image of $f_{0}(\mathfrak{U})$}}
	\end{minipage}
	\hspace{0.1\linewidth}
	\begin{minipage}[h]{0.40\linewidth}
		\centering
		\includegraphics[width=6cm, height=6cm]{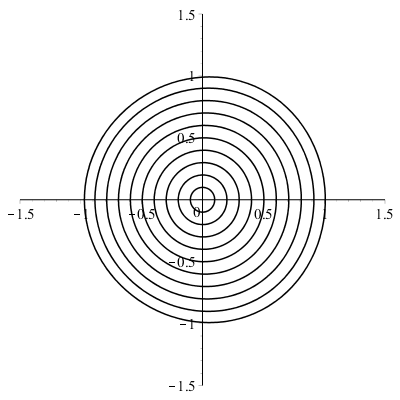}
		\caption{     {\small Image of $P_{0.01,0.01}^{2,2}(f_{0}(\mathfrak{U}))$} }
	
	\end{minipage}
\end{figure}
\end{ex}


\end{document}